\documentclass[12pt]{article}
\usepackage{amsmath,amsthm,amsfonts,amssymb,amscd}
\def\qed{{\unskip\nobreak\hfil\penalty50
\hskip2em\hbox{}\nobreak\hfil$\square$
\parfillskip=0pt \finalhyphendemerits=0\par}\medskip}
\def\proof{\trivlist \item[\hskip \labelsep{\bf Proof\ }]}
\def\endproof{\null\hfill\qed\endtrivlist}

\def\Ad{{\mathrm {Ad}}}

\def\End{{\mathrm {End}}}
\def\Hom{{\mathrm {Hom}}}
\def\id{{\mathrm {id}}}

\def\a{\alpha}

\def\e{\varepsilon}

\def\l{\lambda}

\def\phi{\varphi}

\def\Om{\Omega}

\def\r{{\rho}}

\def\emptyset{\varnothing}
\def\setminus{\smallsetminus}

\def\Diff{{\mathrm {Diff}}}
\def\Mob{{\rm\textsf{M\"ob}}}

\def\Ad{{\mathrm {Ad}}}

\def\End{{\mathrm {End}}}
\def\Hom{{\mathrm {Hom}}}
\def\id{{\mathrm {id}}}

\def\a{\alpha}

\def\e{\varepsilon}

\def\l{\lambda}

\def\phi{\varphi}

\def\Om{\Omega}

\def\r{{\rho}}

\newtheorem{theorem}{Theorem}[section]
\newtheorem{lemma}[theorem]{Lemma}
\newtheorem{conjecture}[theorem]{Conjecture}
\newtheorem{corollary}[theorem]{Corollary}
\newtheorem{definition}[theorem]{Definition}

\newtheorem{proposition}[theorem]{Proposition}
\newtheorem{remark}[theorem]{Remark}

\def\emptyset{\varnothing}
\def\setminus{\smallsetminus}

\def\exp{{\mathrm {exp}}}

\def\Diff{{\mathrm {Diff}}}
\def\Mob{{\rm\textsf{M\"ob}}}

\def\res{\!\restriction\!}
\def\A{{\cal A}}

\def\B{{\cal B}}
\def\C{{\cal C}}

\def\I{{\cal I}}

\def\H{{\cal H}}

\renewcommand{\qed}{\ \hfill $\blacksquare$}

\newcommand{\bdefin}{\begin{definition}}
\newcommand{\blemma}{\begin{lemma}}
\newcommand{\bprop}{\begin{proposition}}
\newcommand{\btheor}{\begin{theorem}}
\newcommand{\bcoro}{\begin{corollary}}
\newcommand{\bconj}{\begin{conjecture}}
\newcommand{\edefin}{\end{definition}}
\newcommand{\elemma}{\end{lemma}}
\newcommand{\eprop}{\end{proposition}}
\newcommand{\etheor}{\end{theorem}}
\newcommand{\ecoro}{\end{corollary}}
\newcommand{\econj}{\end{conjecture}}
\newcommand{\brem}{\begin{remark}}
\newcommand{\erem}{\end{remark}}

\newcommand{\ba}{\begin{array}}
\newcommand{\ea}{\end{array}}
\newcommand{\bea}{\begin{eqnarray}}
\newcommand{\eea}{\end{eqnarray}}
\newcommand{\bean}{\begin{eqnarray*}}
\newcommand{\eean}{\end{eqnarray*}}

\title{\huge Mirror extensions of local nets\\}
\author{
{\sc Feng Xu}\footnote{Supported in part by NSF.}\\
Department of Mathematics\\
University of California at Riverside\\
Riverside, CA 92521\\
E-mail: {\tt xufeng@math.ucr.edu}}
\begin{document}
\date{}
\maketitle

\begin{abstract}
In this paper we prove a general theorem on the extensions of local
nets which was inspired by recent 
examples of exotic extensions for Virasoro nets with central charge
less than one and  earlier work on cosets and conformal inclusions.
When applying the theorem to conformal inclusions and diagonal 
inclusions, we obtain infinite series of new examples of 
completely rational nets.  2000MSC:81R15, 17B69.
\end{abstract}
\newpage

\section{Introduction}
Cosets, orbifolds and simple current extensions are  efficient ways of 
constructing conformal
field theories (CFT), and in fact they are so efficient that many believe
that all rational CFT have been constructed 
(cf. Page 356 of \cite{MS}) \footnote{
Also see Page 256 of \cite{Pol} for a more recent appearance of similar
statement.}.
This also motivates  the following statement of E. Witten from  
Page 356 of \cite{Wi}, that  
``it is possible to conjecture that arbitrary rational conformal
field theories in two dimensions can be derived from Chern-Simon theories
in three dimensions'' .
In this paper we use the
method of operator algebras, especially subfactor theory pioneered by
Vaughan V.F. Jones (cf. \cite{J1}), to produce infinitely many new
rational CFT which do not seem to come from cosets, orbifolds, simple
current extensions or any combinations of them, and they do not 
come from CS theories in the usual way.  
Our results show that one has to modify the above statements
.\par
To describe our main result, 
let $\A$ be a completely rational net (cf. Definition \ref{rational}). 
An irreducible extension of
$\A$ is a net $\B$ such that $\A\subset \B$ is an irreducible subnet
(cf. Definition \ref{ext}). The main result (cf. Th. \ref{main}) in this 
paper 
can be briefly described as follows:
Let $\A\subset \B$ be a normal subnet (cf. Definition \ref{normal}
) and
$\tilde \A$ the coset (cf.  Definition \ref{coset}). Suppose that 
$\A\subset \C$ is an irreducible extension. Then under certain conditions
as described in Th. \ref{main}, there is an irreducible extension
$\tilde \C$ of $\tilde \A$ which is called the ``mirror''  of  
$\A\subset \C$.
\par
That $\tilde \A \subset\tilde \C$ is called the ``mirror'' of
$\A\subset \C$ is explained after Th. \ref{main}. Roughly speaking the
reason is that the link
invariants labeled by the spectrum of $\A\subset \C$ are the invariants
labeled by the spectrum of $\tilde \A \subset\tilde \C$ corresponding to
the mirror image of the link. \par
The motivations for this work come from the papers 
\cite{KL}, \cite{LR}, \cite{Reh2}, \cite{JX}, \cite{Xu1} and \cite{Xu2}. 
In \cite{KL}  two exceptional extensions\footnote{ One of the extensions
is identified as a coset in \cite{Ko}.}
of Virasoro net with central charge $c<1$  were constructed 
based on \cite{LR} 
and known
$A-D-E$ classifications of modular invariants, and in \S4.2 of \cite{Xu1}  
it was observed that two series of conformal inclusions were closely 
related based on the level-rank duality as formulated in \cite{Xu2}.
We will see that both cases are examples of mirror extensions in 
section \ref{examples}. We also note that in  \cite{Reh2} certain
two dimensional local extensions were constructed based on \cite{LR}. 
We will see that the results of \cite{LR} and  
\cite{Xu2} play
an important role in the proof of Th. \ref{main}. The questions at the 
end of \cite{JX} also lead us to develop the general method in this 
paper. \par
Perhaps the simplest example of a mirror extension is the mirror of
conformal inclusion $SU(2)_{10}\subset Spin(5)$. 
This is an extension $\tilde C$ of
$SU(10)_2$, whose spectrum is $L(2\Lambda_0)+ L(\Lambda_3+ \Lambda_7)$ 
(cf. Definition \ref{ext} and subsection \ref{examples} ) 
based on level-rank duality. We note that 
the net  $\tilde C$ is not the net associated to any affine Kac-Moody 
algebras. To the best
of our knowledge, this and most of the infinite series of exotic 
extensions in subsection \ref{examples} have not appeared before, and hence
we obtain a large class of new completely rational net. 
Based on the close relation of nets and  vertex operator algebras (VOA) 
(cf. \cite{Xu3}),  this also
leads to a conjecture about the existence of a large class of 
new rational vertex operator algebras (VOA) (cf. the end of subsection 
\ref{examples}). For
an example, in the case of the simplest example above, the conjecture 
implies the existence of a rational VOA containing the affine VOA based 
on $SU(10)_2$ with  spectrum $L(2\Lambda_0)+ L(\Lambda_3+ \Lambda_7)$. 
We note that this is not an example of 
simple current extensions since the index of the 
inclusion (cf. Definition \ref{ext}) 
is not an integer (it is $3+\sqrt 3$ by \S4.1 of \cite{Xu4}).
\par
The rest of this paper is organized as follows: In \S2 we recall the
basic notions about sectors, nets and examples to set up notations.  
In \S3 we prove our main
result on mirror extensions Th. \ref{main}. Using conformal inclusions,
diagonal embeddings,  and applying  
Th. \ref{main} we obtain a large class of new completely rational nets
in \S4. \par
The author would like to thank 
Professor Vaughan Jones for valuable discussions
and motivation for the problem.

\section{Preliminaries}
\subsection{Preliminaries on sectors}
  
Given an infinite factor $M$, the {\it sectors of $M$}  are given by
$$\text{Sect}(M) = \text{End}(M)/\text{Inn}(M),$$
namely $\text{Sect}(M)$ is the quotient of the semigroup of the
endomorphisms of $M$ modulo the equivalence relation: $\rho,\rho'\in
\text{End}(M),\, \rho\thicksim\rho'$ iff there is a unitary $u\in M$
such that $\rho'(x)=u\rho(x)u^*$ for all $x\in M$.

$\text{Sect}(M)$ is a $^*$-semiring (there are an addition, a product and
an involution $\rho\rightarrow \bar\rho$) 
equivalent to the Connes correspondences (bimodules) on
$M$ up to unitary equivalence. If $\r$ is
an element of $\text{End}(M)$ we shall denote by $[\r]$
its class in $\text{Sect}(M)$. We define
$\text{Hom}(\r,\r')$ between the objects $\r,\r'\in \End(M)$
by 
\[
\text{Hom}(\r,\r')\equiv\{a\in M: a\r(x)=\r'(x)a \ \forall x\in M\}.
\]
We use $\langle  \lambda , \mu \rangle$ to denote
the dimension of   $\text{\rm Hom}(\lambda , \mu )$; it can be 
$\infty$, but it is finite if $\l,\mu$ have finite index. See 
\cite{J1} for the definition of index for type $II_1$ case 
which initiated the subject 
and  \cite{PP} for  the definition of index in general. Also see
\S2.3 
\cite{KLX}
for expositions.  
$\langle  \lambda , \mu \rangle$
depends only on $[\lambda ]$ and $[\mu ]$. Moreover we have
if $\nu$ has finite index, then 
$\langle \nu \lambda , \mu \rangle =
\langle \lambda , \bar \nu \mu \rangle $,
$\langle \lambda\nu , \mu \rangle
= \langle \lambda , \mu \bar \nu \rangle $ which follows from Frobenius
duality.  
$\mu $ is a subsector of $\lambda $ if there is an isometry $v\in M$ 
such that $\mu(x)= v^* \lambda(x)v, \forall x\in M.$
We will also use the following
notation: if $\mu $ is a subsector of $\lambda $, we will write as
$\mu \prec \lambda $  or $\lambda \succ \mu $.  A sector
is said to be irreducible if it has only one subsector. 

\subsection{Local nets}
By an interval of the circle we mean an open connected
non-empty subset $I$ of $S^1$ such that the interior of its 
complement $I'$ is not empty. 
We denote by $\I$ the family of all intervals of $S^1$.

A {\it net} $\A$ of von Neumann algebras on $S^1$ is a map 
\[
I\in\I\to\A(I)\subset B(\H)
\]
from $\I$ to von Neumann algebras on a fixed separable Hilbert space $\H$
that satisfies:
\begin{itemize}
\item[{\bf A.}] {\it Isotony}. If $I_{1}\subset I_{2}$ belong to 
$\I$, then
\begin{equation*}
 \A(I_{1})\subset\A(I_{2}).
\end{equation*}
\end{itemize}
If $E\subset S^1$ is any region, we shall put 
$\A(E)\equiv\bigvee_{E\supset I\in\I}\A(I)$ with $\A(E)=\mathbb C$ 
if $E$ has empty interior (the symbol $\vee$ denotes the von Neumann 
algebra generated). 

The net $\A$ is called {\it local} if it satisfies:
\begin{itemize}
\item[{\bf B.}] {\it Locality}. If $I_{1},I_{2}\in\I$ and $I_1\cap 
I_2=\emptyset$ then 
\begin{equation*}
 [\A(I_{1}),\A(I_{2})]=\{0\},
 \end{equation*}
where brackets denote the commutator.
\end{itemize}
The net $\A$ is called {\it M\"{o}bius covariant} if in addition 
satisfies
the following properties {\bf C,D,E,F}:
\begin{itemize}
\item[{\bf C.}] {\it M\"{o}bius covariance}. 
There exists a non-trivial strongly 
continuous unitary representation $U$ of the M\"{o}bius group 
$\Mob$ (isomorphic to $PSU(1,1)$) on $\H$ such that
\begin{equation*}
 U(g)\A(I) U(g)^*\ =\ \A(gI),\quad g\in \Mob,\ I\in\I.
\end{equation*}
\item[{\bf D.}] {\it Positivity of the energy}.  
The generator of the one-parameter
rotation subgroup of $U$ (conformal Hamiltonian), denoted by 
$L_0$ in the following,  is positive.
\item[{\bf E.}] {\it Existence of the vacuum}.  There exists a unit
$U$-invariant vector $\Omega\in\H$ (vacuum vector), and $\Omega$ is
cyclic for the von Neumann algebra $\bigvee_{I\in\I}\A(I)$.
\end{itemize}
By the Reeh-Schlieder theorem $\Omega$ is cyclic and separating for 
every fixed $\A(I)$. The modular objects associated with 
$(\A(I),\Omega)$ have a geometric meaning
\[
\Delta^{it}_I = U(\Lambda_I(2\pi t)),\qquad J_I = U(r_I)\ .
\]
Here $\Lambda_I$ is a canonical one-parameter subgroup of $\Mob$ and $U(r_I)$ is a 
antiunitary acting geometrically on $\A$ as a reflection $r_I$ on $S^1$. 

This implies {\em Haag duality}: 
\[
\A(I)'=\A(I'),\quad I\in\I\ ,
\]
where $I'$ is the interior of $S^1\setminus I$.

\begin{itemize}
\item[{\bf F.}] {\it Irreducibility}. $\bigvee_{I\in\I}\A(I)=B(\H)$. 
Indeed $\A$ is irreducible iff
$\Om$ is the unique $U$-invariant vector (up to scalar multiples). 
Also  $\A$ is irreducible
iff the local von Neumann 
algebras $\A(I)$ are factors. In this case they are either ${\mathbb C}$ or 
III$_1$-factors 
with separable predual
in 
Connes classification of type III factors.
\end{itemize}
By a {\it conformal net} (or diffeomorphism covariant net)  
$\A$ we shall mean a M\"{o}bius covariant net such that the following 
holds:
\begin{itemize}
\item[{\bf G.}] {\it Conformal covariance}. There exists a projective 
unitary representation $U$ of $\Diff(S^1)$ on $\H$ extending the unitary 
representation of $\Mob$ such that for all $I\in\I$ we have
\begin{gather*}
 U(\phi)\A(I) U(\phi)^*\ =\ \A(\phi.I),\quad  \phi\in\Diff(S^1), \\
 U(\phi)xU(\phi)^*\ =\ x,\quad x\in\A(I),\ \phi\in\Diff(I'),
\end{gather*}
\end{itemize}
where $\Diff(S^1)$ denotes the group of smooth, positively oriented 
diffeomorphism of $S^1$ and $\Diff(I)$ the subgroup of 
diffeomorphisms $g$ such that $\phi(z)=z$ for all $z\in I'$.
\par
A (DHR) representation $\pi$ of $\A$ on a Hilbert space $\H$ is a map 
$I\in\I\mapsto  \pi_I$ that associates to each $I$ a normal 
representation of $\A(I)$ on $B(\H)$ such that
\[
\pi_{\tilde I}\res\A(I)=\pi_I,\quad I\subset\tilde I, \quad 
I,\tilde I\subset\I\ .
\]
$\pi$ is said to be M\"obius (resp. diffeomorphism) covariant if 
there is a projective unitary representation $U_{\pi}$ of $\Mob$ (resp. 
$\Diff(S^1)$) on $\H$ such that
\[
\pi_{gI}(U(g)xU(g)^*) =U_{\pi}(g)\pi_{I}(x)U_{\pi}(g)^*
\]
for all $I\in\I$, $x\in\A(I)$ and $g\in \Mob$ (resp. 
$g\in\Diff(S^1)$). 

By definition the irreducible conformal net is in fact an irreducible 
representation of itself and we will call this representation the {\it 
vacuum representation}.\par

Let $G$ be a simply connected  compact Lie group. By Th. 3.2
of \cite{FG}, 
the vacuum positive energy representation of the loop group
$LG$ (cf. \cite{PS}) at level $k$ 
gives rise to an irreducible conformal net 
denoted by {\it ${\A}_{G_k}$}. By Th. 3.3 of \cite{FG}, every 
irreducible positive energy representation of the loop group
$LG$ at level $k$ gives rise to  an irreducible covariant representation 
of ${\A}_{G_k}$. \par
Given an interval $I$ and a representation $\pi$ 
of $\A$, there is an {\em endomorphism of $\A$ localized in $I$} equivalent 
to $\pi$; namely $\r$ is a representation of $\A$ on the vacuum Hilbert 
space $\H$, unitarily equivalent 
to $\pi$, such that $\r_{I'}=\text{id}\restriction\A(I')$.
We now define  the statistics. Given the endomorphism $\r$ of $\A$
localized in $I\in\I$, choose an equivalent endomorphism $\r_0$
localized in an interval $I_0\in\I$ with $\bar I_0\cap\bar I
=\emptyset$ and let $u$ be a local intertwiner in $\Hom(\r,\r_0)$ , 
namely $u\in \Hom(\r_{\tilde I},\r_{0,\tilde I})$ with $I_0$ following
clockwise $I$ inside $\tilde I$ which is an interval containing
both $I$ and $I_0$.

The {\it statistics operator} $\epsilon (\r,\rho):= u^*\r(u) =
u^*\r_{\tilde I}(u) $ belongs to $\Hom(\r^2_{\tilde I},\r^2_{\tilde
I})$. We will call $\epsilon (\r,\rho)$ the positive or right braiding and
$\tilde\epsilon (\r,\rho):=\epsilon (\r,\rho)^*$ the negative or left 
braiding.
\par  
Next we  recall some definitions from \cite{KLM} .
Recall that   ${\I}$ denotes the set of intervals of $S^1$.
Let $I_1, I_2\in {\I}$. We say that $I_1, I_2$ are disjoint if
$\bar I_1\cap \bar I_2=\emptyset$, where $\bar I$
is the closure of $I$ in $S^1$.  
When $I_1, I_2$ are disjoint, $I_1\cup I_2$
is called a 1-disconnected interval in \cite{Xu3}.  
Denote by ${\I}_2$ the set of unions of disjoint 2 elements
in ${\I}$. Let ${\A}$ be an irreducible M\"{o}bius covariant net
. For $E=I_1\cup I_2\in{\I}_2$, let
$I_3\cup I_4$ be the interior of the complement of $I_1\cup I_2$ in 
$S^1$ where $I_3, I_4$ are disjoint intervals. 
Let 
$$
{\A}(E):= A(I_1)\vee A(I_2), \quad
\hat {\A}(E):= (A(I_3)\vee A(I_4))'.
$$ Note that ${\A}(E) \subset \hat {\A}(E)$.
Recall that a net ${\A}$ is {\it split} if ${\A}(I_1)\vee
{\A}(I_2)$ is naturally isomorphic to the tensor product of
von Neumann algebras ${\A}(I_1)\otimes
{\A}(I_2)$ for any disjoint intervals $I_1, I_2\in {\I}$.
${\A}$ is {\it strongly additive} if ${\A}(I_1)\vee
{\A}(I_2)= {\A}(I)$ where $I_1\cup I_2$ is obtained
by removing an interior point from $I$.
\bdefin\label{rational}
\cite{{KLM}, {LX}}
A M\"{o}bius covariant net ${\A}$ is said to be completely  rational if
${\A}$ is split, and 
the index $[\hat {\A}(E): {\A}(E)]$ is finite for some
$E\in {\I}_2$ . The value of the index
$[\hat {\A}(E): {\A}(E)]$ (it is independent of 
$E$ by Prop. 5 of \cite{KLM}) is denoted by $\mu_{{\A}}$
and is called the $\mu$-index of ${\A}$. 
\edefin
Note that, by recent results in \cite{LX}, 
every irreducible, split, 
local conformal net with finite $\mu$-index is automatically strongly 
additive. Hence we have modified the definition in \cite{KLM} by dropping the
strong additivity requirement in the above definition. Also note that
if $\A$ is completely rational, then $\A$ has only finitely many 
irreducible covariant representations by \cite{KLM}.
\par
Let $\B$ be a  M\"{o}bius net. By a {\it  M\"{o}bius subnet} (cf. \cite{L1}) 
we shall mean a
map  
\[
I\in\I\to\A(I)\subset \B(I)
\]
that associates to each interval $I\in \I$ a von Neumann subalgebra $\A(I)$
of $\B(I)$, which is isotonic
\[
\A(I_1)\subset \A(I_2), I_1\subset I_2,
\]
and   M\"{o}bius covariant with respect to the the representation $U$, 
namely 
\[
U(g) \A(I) U(g)^*= \A(g.I)
\] for all $g\in \Mob$ and $I\in \I$. Note that by Lemma 13 
of \cite{L1} for each $I\in \I$ there exists a conditional 
expectation $E_I: \B(I)\rightarrow \A(I)$ such that $E$ preserves the
vector state given by the vacuum of $\A$.
\begin{definition}\label{ext}
Let $\A$ be a  M\"{o}bius covariant net. A  M\"{o}bius covariant net $\B$
on a Hilbert space $\H$ is an extension of $\A$ if there is a DHR
representation  $\pi$ of $\A$ on $\H$ such that $\pi(\A)\subset \B$
is a   M\"{o}bius subnet. The extension is irreducible if 
$\pi(\A(I))'\cap \B(I) = {\mathbb C} $ for some (and hence all) interval $I$,
and is of finite index if $\pi(\A(I))\subset \B(I)$ has finite 
index for some (and hence all) interval $I$. The index will be called the
index of the inclusion $\pi(\A)\subset \B$. If $\pi$ as representation of
 $\A$ decomposes as $[\pi]= \sum_\lambda m_\lambda[\lambda]$ where 
$m_\lambda$  are non-negative  integers and $\lambda$ are irreducible
DHR representations of $\A$, we say that 
$[\pi]= \sum_\lambda m_\lambda[\lambda]$ is the spectrum of the
extension.  For simplicity we will write  $\pi(\A)\subset \B$ simply as
$\A\subset \B$.
\end{definition}
\begin{lemma}\label{ext1}
If  $\A$ is completely rational, and a   M\"{o}bius covariant net $\B$ is 
an irreducible extension of $\A$. Then $\A\subset\B$ has finite 
index and $\B$ is completely rational. 
\end{lemma}
\proof
 $\A\subset\B$ has finite 
index follows from Prop. 2.3 of \cite{KL}, and it follows by 
Prop. 24 of \cite{KLM} that   $\B$ is completely rational. 
\endproof
The following is essentially Th. 4.9 of \cite{LR} (cf. \S2.4 of \cite{KL})
which is also used in \S4.2 of  \cite{KL}:
\begin{proposition}\label{qlocal}
Let $\A$ be a  M\"{o}bius covariant net, $\r$ a DHR representation 
of $\A$ localized on a fixed $I_0$ with finite statistics, which contains 
$\id$ with multiplicity one, i.e., there is (unique up to a phase) 
isometry $w\in\Hom(\id,\r).$ Then there is a  M\"{o}bius 
covariant net $\B$ which is an irreducible extension of $\A$ if and only
if there is an isometry $w_1\in\Hom(\r,\r^2)$ 
which solves the following equations:
\begin{align}
w_1^* w & = w_1^* \r(w) \in {\mathbb R_+} \label{a}\\
w_1w_1& = \r(w_1) w_1 \label{b} \\
\epsilon(\r,\r) w_1 & = w_1 \label{c}
\end{align}
\end{proposition}
\proof
As in the proof of Th. 4.9 in \cite{LR}, we just have to check the
``if'' part, and the only additional thing we need to check is that
$\B$ is  M\"{o}bius  covariant. Since  $\r$ is DHR with finite
statistics, it follows that  $\r$ is  M\"{o}bius  covariant by 
\cite{GL}, and it follows by the formula in Cor. 19 of \cite{L1} that
$\B$ is  M\"{o}bius  covariant.
\endproof
\subsection{Induction}\label{ind}
Let $\B$ be a M\"obius covariant net and $\A$ a subnet. 
We assume that $\A$ is strongly additive and $\A\subset \B$ has
finite index.  
Fix an 
interval $I_0\in\I$ and  canonical endomorphism (cf. \cite{LR})
$\gamma$ associated with $\A(I_0)\subset\B(I_0)$. 
Given a DHR endomorphism $\r$ of $\B$ localized in $I_0$, the 
$\a$-induction $\a_{\r}$ 
of $\r$ is the endomorphism of $\B(I_0)$ given by
\[
\a_{\r}\equiv \gamma^{-1}\cdot\Ad\e(\r,\l)\cdot\r\cdot\gamma\ 
\]
where $\e$ denotes the right braiding  (cf. Cor. 3.2 of \cite{BE}).
In \cite{Xu4} a slightly different endomorphism was introduced and the
relation between the two was given in \S2.1 of \cite{Xu2}.
Note that $\Hom( \a_\lambda,\a_\mu)=:\{ x\in \B(I_0) |
x  \a_\lambda(y)= \a_\mu(y)x, \forall y\in \B(I_0)\} $ and
$\Hom( \lambda,\mu)= :\{ x\in \A(I_0) |
x  \lambda(y)= \mu(y)x, \forall y\in \A(I_0)\} .$
We have the following lemma which will be used in the proof of
Prop. \ref{mirror}. The proof is also implicitly contained in 
 \cite{Xu4}.
\begin{lemma}\label{A=B}  
Assume that $\lambda,\mu$ have finite index  and 
$\langle \lambda, \mu\rangle= 
\langle \a_\lambda, \a_\mu\rangle$. Then 
$$
\Hom( \a_\lambda,\a_\mu)=\Hom( \lambda,\mu)\subset  \A(I_0).$$
\end{lemma}
\proof
By Lemma 3.6 of \cite{BE} $\Hom( \lambda,\mu)\subset
\Hom( \a_\lambda,\a_\mu)$, and by assumption
$$
\langle \lambda, \mu\rangle= 
\langle \a_\lambda, \a_\mu\rangle < \infty,
$$ 
hence
$$\Hom( \a_\lambda,\a_\mu)=\Hom( \lambda,\mu)\subset  \A(I_0).$$
\endproof
\section{Mirror extensions}
\subsection{Coset construction}
Let $\B$ be a completely rational net and $\A\subset \B$ be a subnet 
which is also completely rational. 
\begin{definition}\label{coset}
Define a subnet $\tilde\A\subset \B$ by 
$\tilde\A(I):= \A(I)'\cap \B(I), \forall I\in \I.$
\end{definition}
We note that since $\A$ is completely rational, it is strongly additive
and so we have 
$\tilde\A(I)= (\vee_{J\in \I}\A(J))'\cap \B(I), \forall I\in \I.$
The following lemma then follows directly from the definition:
\begin{lemma}\label{cosetnet}
The restriction of $\tilde\A$ on the Hilbert space $
\overline{ \vee_I\tilde\A(I)\Omega}$ is an irreducible
  M\"{o}bius  covariant net.
\end{lemma}
The net $\tilde\A$ as in Lemma \ref{cosetnet} will be called
the {\it coset} of  $\A\subset \B$. See \cite{Xu3} for 
a class of cosets from Loop groups. \par
The following definition generalizes the definition in \S3 of
\cite{Xu3}:
\begin{definition}\label{cofinite}  
 $\A\subset \B$ is called cofinite if the inclusion
$\tilde\A(I)\vee \A(I) \subset \B(I)$ has finite 
index for some interval $I$.
\end{definition}
\begin{proposition}\label{rationalc}
Let $\B$ be completely rational, and 
let $\A\subset\B$ be a  M\"{o}bius subnet which is also completely
rational. Then $\A\subset \B$ is  cofinite if and only if
$\tilde\A$ is completely rational. 
\end{proposition}
\proof
Note that by \cite{SS} $\tilde\A(I)\vee \A(I)$ is naturally 
isomorphic to $\tilde \A(I)\otimes \A(I)$, and by the geometric nature of the
modular group $\tilde\A(I)$ is 
a factor. It follows that the inclusion $\tilde\A(I)\vee \A(I) \subset \B(I)$
is irreducible. The ``if'' part now follows from Prop. 2.3 of \cite{KL},
and the ``only if'' part follows from Prop. 24 of \cite{KLM}.
\endproof
Let  $\B$ be completely rational, and 
let $\A\subset\B$ be a  M\"{o}bius subnet which is also completely
rational. Assume that  $\A\subset \B$ is cofinite. 
We will use $\sigma_i,\sigma_j,...$ (resp. $\lambda, \mu...$) 
to label irreducible DHR representations of 
$\B$ (resp. $\A$) localized on a fixed interval $I_0$. 
Since $ \tilde\A$ is completely rational by 
Prop. \ref{rationalc}, $\tilde\A\otimes \A$ is completely rational, 
and so every irreducible DHR representation $\sigma_i$ of $\B$, when
restricting to $\tilde\A\otimes \A$, decomposes as direct sum of 
representations of  $ \tilde\A\otimes \A$ of the form 
$(i,\lambda)\otimes \lambda$ by Lemma 27 of \cite{KLM}. Here $(i,\lambda)$
is a DHR representation of $\tilde \A$ which may not be irreducible and we 
use the  tensor notation $(i,\lambda)\otimes \lambda$ to represent a 
DHR representation of $ \tilde\A\otimes \A$ which is localized on
$I_0$ and defined by
$$
(i,\lambda)\otimes \lambda (x_1\otimes x_2)= (i,\lambda)(x_1)\otimes 
 \lambda(x_2), \forall x_1\otimes x_2\in \tilde\A(I_0) \otimes \A(I_0).
$$
We will also identify $\tilde \A$ and $\A$ as subnets of 
$ \tilde\A\otimes \A$ in the natural way. 
We note that when no confusion arise, we will use $1$ to denote the
vacuum representation of a net. 
  
\begin{definition}\label{normal}
A  M\"{o}bius subnet $\A\subset\B$ is normal if 
$\tilde \A(I)'\cap \B(I)= \A$ for some I.
\end{definition}
The following is implied by Lemma 3.4 of \cite{Reh1} (also cf. Page 797
of \cite{Xu5}):
\begin{lemma}\label{normal1}
Let  $\B$ be completely rational, and 
let $\A\subset\B$ be a  M\"{o}bius subnet which is also completely
rational. Assume that  $\A\subset \B$ is cofinite. Then the following 
conditions are equivalent: \par
(1) $\A\subset\B$ is normal; \par
(2) $(1,1)$ is the vacuum representation of $\tilde \A$ and 
$(1,\lambda) $ contains $(1,1)$ if and only if $\lambda=1$. \par
\end{lemma}
The following proposition will play a key role in the proof of
Th. \ref{main}:
\begin{proposition}\label{mirror}
Let  $\B$ be completely rational, and 
let $\A\subset\B$ be a  M\"{o}bius subnet which is also completely
rational. Assume that  $\A\subset \B$ is cofinite and normal. Then:\par
(1) Let $\gamma$ be the restriction of the vacuum representation of $\B$
to $\tilde\A\otimes \A$. 
Then $[\gamma]= \sum_{\lambda\in \exp} [(1,\lambda)\otimes \lambda]$
where each $(1,\lambda)$ is irreducible;\par
(2) Let $\lambda\in \exp$ be as in (1), then $[\alpha_{(1,\lambda)\otimes 1}]
= [\alpha_{1\otimes \bar\lambda}]$, and $[\lambda]\rightarrow 
[\alpha_{1\otimes \lambda}]$ is a ring isomorphism where the $\a$-induction
is with respect to $\tilde\A\otimes \A \subset \B$ as in subsection 
\ref{ind}; Moreover   
the set $\exp$ is closed under fusion;
\par
(3) Let $[\rho]= \sum_{\lambda\in \exp} m_\lambda[\lambda]$ where 
$m_\lambda = m_{\bar\lambda}\geq 0, \forall \lambda$, and $
[(1,\rho)]= \sum_{\lambda\in \exp} m_\lambda[(1,\lambda)]$. Then there
exists an unitary element $T_\rho \in
\Hom(\alpha_{(1,\rho)\otimes 1}, 
\alpha_{1\otimes \rho})$ such that
$$
\epsilon((1,\rho),(1,\rho)) T_\rho^* \alpha_{1\otimes \rho}(T_\rho^*)
=T_\rho^* \alpha_{1\otimes \rho}(T_\rho^*) \tilde\epsilon(\rho,\rho)
;$$ \par
(4) Let $\rho$, $(1,\rho)$ be as in (3). Then 
\begin{align*}
\Hom(\rho^n, \rho^m) & = \Hom(\alpha_{1\otimes \rho^n},\alpha_{1\otimes \rho^m}),  \\
\Hom((1,\rho)^n, (1,\rho)^m) & = \Hom(\alpha_{(1,\rho)^n\otimes 1},
\alpha_{(1,\rho)^m\otimes 1}), \forall n,m\in {\mathbb N};
\end{align*}
\par

(5)
Let $\rho$, $(1,\rho)$ be as in (3), let $L$ be an oriented framed
link in three sphere with $n$ components and let $\tilde L$ be the
mirror image of $L$ (cf. \cite{Tu}). Then 
$\tilde L(\rho,...,\rho)= L((1,\rho),...,(1,\rho))= 
\overline{ L(\rho,...,\rho)} $, where $L(i_1,...i_k)$ is defined
as before Lemma 1.7.4 of \cite{Xu2}.
\end{proposition}
\proof
(1),(2) follow from Th. 3.6 of \cite{Reh1} (also cf. Prop. 4.3 of \cite{Xu3}).
As for (3), note that $[\bar\rho]=[\rho]$, and by (2) there exists
an unitary element $T_\rho: \alpha_{(1,\rho)\otimes 1}\rightarrow 
\alpha_{1\otimes \rho}.$  The  equation in (3) follows from 
(3) of Prop. 2.3.1 in \cite{Xu2}. (4) follows from 
(2) and Lemma \ref{A=B}.
(5) follows from Th. B in \cite{Xu2}: we note that even though 
Th. B of  \cite{Xu2} is stated for cosets coming from Loop groups,
the proof of Th. B applies verbatim to our case.
\endproof
\begin{theorem}\label{main}
Let  $\B$ be completely rational, and 
let $\A\subset\B$ be a  M\"{o}bius subnet which is also completely
rational. Assume that  $\A\subset \B$ is cofinite and normal, and let $\exp$
be as in (1) of Prop.\ref{mirror}.  Assume that 
$\A\subset\C$ is an irreducible   M\"{o}bius extension of $\A$ with 
spectrum  $[\rho]=\sum_{\lambda\in \exp} 
m_\lambda [\lambda], m_\lambda\geq 0.$ 
Then there is  an irreducible   M\"{o}bius extension  
$\tilde \C$ of $\tilde\A$ with spectrum 
$[(1,\rho)]=\sum_{\lambda\in \exp} m_\lambda [(1,\lambda)]$.
Moreover $\tilde \C$ is completely rational.
\end{theorem}
\proof
Since $\A\subset\C$ is an irreducible   M\"{o}bius extension of $\A$ with 
spectrum  $[\rho]=\sum_{\lambda\in \exp} 
m_\lambda [\lambda], m_\lambda\geq 0$, by Prop. \ref{qlocal} there
exist $w\in \Hom(\id , \rho), w_1\in \Hom(\rho,\rho^2)$ which
verifies equations (\ref{a}),(\ref{b}) and (\ref{c}). 
Note that $[\rho]=[\bar\rho]$ by \cite{LR}. 
Let $T_\rho$
be the unitary as given by (3) of Prop. \ref{mirror} and define
$\tilde w:= T_\rho^* w, \tilde w_1:= T_\rho^* 
\alpha_{1\otimes \rho}(T_\rho^*) w_1 T_\rho$. Note that by definitions
$$\tilde w\in \Hom(1,\alpha_{(1,\rho)\otimes 1})= \Hom(1, (1,\rho)), 
 \tilde w_1\in \Hom(\alpha_{(1,\rho)\otimes 1}, \alpha_{(1,\rho)\otimes 1}^2)
= \Hom((1,\rho), (1,\rho)^2)
$$
where we have also used (4) of Prop. \ref{mirror}. To prove the theorem,
by Prop. \ref{mirror} and Lemma \ref{ext1} 
it is enough to check equations (\ref{a}), (\ref{b}) and (\ref{c}) with 
$\r$ replaced by $(1,\rho)$ in Prop. \ref{qlocal}. First let us check
 equation (\ref{a}): 
$$
\tilde w_1^* \tilde w=  T_\rho^* w_1^*\alpha_{1\otimes \rho}(T_\rho) 
T_\rho T_\rho^* w
= T_\rho^* w_1^*\alpha_{1\otimes \rho} (T_\r)w
=  T_\rho^*  w_1^* w T_\rho
=  w_1^* w
$$

where we have used $w\in \Hom(1,\rho)= \Hom(1,\alpha_{1\otimes \rho})$ and 
$ w_1^* w\in {\mathbb R_+}$. similarly 
$$
\tilde w_1^* (1,\rho) (\tilde w)=  w_1^* \rho(w)\in {\mathbb R_+}.
$$ 
Next we have
\begin{align*}
\tilde w_1 \tilde w_1 &= T_\rho^*\alpha_{1\otimes \rho}(T_\rho^*) w_1 T_\rho
 T_\rho^*\alpha_{1\otimes \rho}(T_\rho^*) w_1 T_\rho\\
&=  T_\rho^*\alpha_{1\otimes \rho}(T_\rho^*) w_1 
\alpha_{1\otimes \rho}(T_\rho^*) w_1 T_\rho\\
&= T_\rho^*\alpha_{1\otimes \rho}(T_\rho^*) \alpha_{1\otimes \rho}^2(T_\rho^*)
w_1^2 T_\rho\\
&=T_\rho^*\alpha_{1\otimes \rho}(T_\rho^*) \alpha_{1\otimes \rho}^2(T_\rho^*)
\alpha_{1\otimes \rho}(w_1)w_1  T_\rho
\end{align*}
where we have used $w_1\in \Hom(\alpha_{1\otimes \rho}, \alpha_{1\otimes \rho}^2)
=\Hom(\rho,\rho^2), w_1^2= \rho(w_1)w_1$ and  $\alpha_{1\otimes \rho}(w_1)=
\rho(w_1)$ by definition. On the other hand
\begin{align*}
(1,\rho)(\tilde w_1) \tilde w_1 &=
\alpha_{(1,\rho)\otimes 1}( T_\rho^*\alpha_{1\otimes \rho}(T_\rho^*) w_1 T_\rho) T_\rho^*\alpha_{1\otimes \rho}(T_\rho^*) w_1 T_\rho \\
& = \alpha_{(1,\rho)\otimes 1}( T_\rho^*\alpha_{1\otimes \rho}(T_\rho^*))
\alpha_{(1,\rho)\otimes 1}(w_1)  \alpha_{(1,\rho)\otimes 1}(T_\rho)
 T_\rho^*\alpha_{1\otimes \rho}(T_\rho^*) w_1 T_\rho \\
& = \alpha_{(1,\rho)\otimes 1}( T_\rho^*\alpha_{1\otimes \rho}(T_\rho^*))
\alpha_{(1,\rho)\otimes 1}(w_1) T_\rho^* w_1 T_\rho \\
& = \alpha_{(1,\rho)\otimes 1}( T_\rho^*\alpha_{1\otimes \rho}(T_\rho^*))
 T_\rho^*\alpha_{1\otimes \rho}(w_1) w_1 T_\rho \\
&=T_\rho^*\alpha_{1\otimes \rho}(T_\rho^*) \alpha_{1\otimes \rho}^2(T_\rho^*)
\alpha_{1\otimes \rho}(w_1)w_1  T_\rho
\end{align*}
This proves   equation (\ref{b}).  
By (3) of Prop. \ref{mirror}, we have
$$
\epsilon((1,\rho),(1,\rho)) \tilde w_1= 
\epsilon((1,\rho),(1,\rho))  T_\rho^*\alpha_{1\otimes \rho}(T_\rho^*)
w_1 T_\rho =
 T_\rho^*\alpha_{1\otimes \rho}(T_\rho^*) \tilde\epsilon(\rho,\rho) 
w_1 T_\rho,
$$
and since $\epsilon(\rho,\rho)w_1=w_1, \tilde\epsilon(\rho,\rho)=
\epsilon(\rho,\rho)^*$, we have proved 
$\epsilon((1,\rho),(1,\rho)) \tilde w_1=  \tilde w_1$ 
which is equation (\ref{c}).
\endproof
\begin{remark}\label{mc}
Due to (5) of Prop. \ref{mirror}, the extension $\tilde \A\subset\tilde \C$
as given in Th. \ref{main} will be called the mirror or the conjugate of
$\A\subset \C$.
\end{remark} 
\begin{remark}\label{nonlocal}
The same idea in the proof of  Th. \ref{main} can also be used to obtain 
possibly non-local  extension $\tilde \A\subset\tilde \C$ when 
$\B$ is not necessarily local, and we plan to discuss applications
elsewhere.
\end{remark} 
\section{Applications}
\subsection{Extensions from conformal inclusions}
Let $G= SU(N)$. We denote $LG$ the group of smooth maps
$f: S^1 \mapsto G$ under pointwise multiplication. The
diffeomorphism group of the circle $\text{\rm Diff} S^1 $ is 
naturally a subgroup of $\text{\rm Aut}(LG)$ with the action given by 
reparametrization. In particular the group of rotations
$\text{\rm Rot}S^1 \simeq U(1)$ acts on $LG$. We will be interested 
in the projective unitary representation $\pi : LG \rightarrow U(H)$ that 
are both irreducible and have positive energy. This means that $\pi $ 
should extend to $LG\ltimes \text{\rm Rot}\ S^1$ so that
$H=\oplus _{n\geq 0} H(n)$, where the $H(n)$ are the eigenspace
for the action of $\text{\rm Rot}S^1$, i.e.,
$r_\theta \xi = \exp^{i n \theta}$ for $\theta \in H(n)$ and 
$\text{\rm dim}\ H(n) < \infty $ with $H(0) \neq 0$. It follows from 
\cite{PS} that for fixed level $K$ which
is a positive integer, there are only finite number of such 
irreducible representations indexed by the finite set
$$
 P_{++}^{K} 
= \bigg \{ \lambda \in P \mid \lambda 
= \sum _{i=1, \cdots , N-1}
\lambda _i \Lambda _i , \lambda _i \geq 0\, ,
\sum _{i=1, \cdots , n-1}
\lambda _i \leq K \bigg \}
$$
where $P$ is the weight lattice of $SU(N)$ and $\Lambda _i$ are the 
fundamental weights. We will use 
$\Lambda_0$ to denote the trivial representation of 
$SU(N)$. We will use $L(\lambda),  \lambda\in P_{++}^{K}$ to label the 
irreducible representations of the net $\A_{SU(N)_K}$.\par  
Let $G \subset H$ be inclusions of compact simple Lie groups.  $LG 
\subset LH$ is called a conformal inclusion if the level 1 projective positive 
energy representations of $LH$ decompose as a finite number of irreducible 
projective representations of $LG$.   $LG
\subset LH$ is called a maximal  conformal inclusion
if there is no proper subgroup $G'$ of $H$ containing 
$G$ such that   $LG
\subset LG'$ is also  a conformal inclusion.
A list of maximal conformal inclusions can be 
found in \cite{GNO}.  \par
Let $H^0$ be the vacuum representation of $LH$, i.e., the
 representation of $LH$ associated with the trivial  representation of $H$.  
Then $H^0$  
decomposes as a direct sum of  irreducible projective representation  
of $LG$ at level $K$. $K$  is called the Dynkin index of the 
conformal inclusion.

We shall  write the 
conformal inclusion as $G_K\subset H_1$. Note that it follows 
from the definition that $\A_{H_1}$ is an extension of
$\A_{G_K}$.
We shall limit our consideration to the following conformal 
inclusions so we can 
use the results of \cite{Xu3}, \cite{Xu4} which is 
based on \cite{W}, though most of the arguments apply to other cases as well
under certain finite index assumptions.\par
\noindent
$\ \ {SU}(2)_{10} \subset {SO}(5)_1, \ 
{SU}(2)_{28} \subset (G_2)_1$,
 
\noindent
$\ \ {SU}(3)_5 \subset {SU}(6)_1, \ 
{SU}(3)_9 \subset (E_6)_1, \ {SU}(3)_{21} \subset 
(E_7)_1$; \par
\noindent
$(A_8)_1 \subset (E_8)_1$;\par
\noindent
and four infinite series:
\begin{align}
{SU}(N)_{N-2} & \subset \ {SU} \left( \frac{N(N-1)}{2} 
\right)_1, \ \ N \geq 4 ; \label{c1} \\
{SU}(N)_{N+2} & \subset \ {SU} \left( \frac{N(N+1)}{2} 
\right)_1;  \label{c2} \\
{SU}(N)_{N} & \subset \ {Spin} (N^2 - 1)_1, \ \ N \geq 
2;  \label{c3} \\
{SU}(N)_M \times {SU}(M)_N  & \subset \ {SU(MN)}_1.
\label{c4} 
\end{align}
Note that except equation (\ref{c4}), the above  
cover all the maximal conformal inclusions of the 
form $SU(N) \subset H$ with
$H$ being a simple group. 
\subsection{Two series of normal inclusions}
\begin{lemma}\label{lr}
The subnets $\A_{{SU}(N)_M}\subset \A_{{SU}(NM)_1}$
are normal and cofinite.  the set $\exp$ as in (1) Prop. \ref{mirror} is
the elements of $P_{++}^{N+M}$ which belong to the root lattice of 
${SU}(N)$.
\end{lemma}
\proof
By lemma 3.3 of \cite{Xu3} the coset of  
$\A_{{SU}(N)_M}\subset \A_{{SU}(NM)_1}$ 
can be identified with $\A_{{SU}(M)_N}$, and exchanging $M$ and $N$ we conclude that
 $\A_{{SU}(N)_M}\subset \A_{{SU}(NM)_1}$ is normal, and it 
is cofinite by the remark after lemma 3.3 in \cite{Xu3}. The statement about 
$\exp$ follows from the branching rules in \cite{ABI}.
\endproof
\begin{lemma}\label{diagonal}
(1) If $\A\subset\B$ and $\B\subset \C$ are normal subnets, then
$\A\subset \C$ is also normal; \par
(2) Let $\A_{{SU}(N)_{K}} \subset \A_{{SU}(N)_{K_1}}\otimes
\A_{{SU}(N)_{K_2}}\otimes...\otimes \A_{{SU}(N)_{K_l}}$
be subnets corresponding to the diagonal embedding of 
${SU}(N)$ in ${SU}(N)\times...\times {SU}(N)$
($l$ tensor factors), where $K=K_1+...+K_l$. Then 
 $\A_{{SU}(N)_{K}} \subset \A_{{SU}(N)_{K_1}}\otimes
\A_{{SU}(N)_{K_2}}\otimes...\otimes \A_{{SU}(N)_{K_l}}$ 
is normal
and cofinite. Moreover, the set $\exp$ as in (1) Prop. \ref{mirror} is
the elements of $P_{++}^K$ which belong to the root lattice of 
${SU}(N)$.
\end{lemma}
\proof
Ad (1):
For a fixed interval $I$, we have
$\A(I)'\cap \C(I)\supset \A(I)'\cap \B(I) \vee \B(I)'\cap \C(I)$, hence
$$
(\A(I)'\cap \C(I))' \cap \C(I)\subset 
(\A(I)'\cap \B(I))'\cap  (\B(I)'\cap \C(I))'\cap \C(I)
=(\A(I)'\cap \B(I))'\cap \B(I)
=\A(I)
$$ by the normality of  $\A\subset\B$ and $\B\subset \C$. It follows that
$ (\A(I)'\cap \C(I))' \cap \C(I) =\A(I)$ and $\A\subset \C$ is normal. \par
Ad(2): By Cor. 3.4 of \cite{Xu3} the subnet is cofinite, so we just have to 
show the normality. 
We first show (2) for the case $l=2$. By lemma \ref{normal1} it
is sufficient to show that a representation $(1,\lambda)$ contains the 
vacuum representation $(1,1)$ if and only if $\lambda=1$, and this 
follows from the remark on Page 38 of \cite{Xu3} before equation (**). 
The statement about $\exp$ follows by Page 194 of \cite{KW}. The general 
$l$ case follows by induction from the case $l=2$ and (1).
\endproof
\subsection{Mirror extensions corresponding to  conformal inclusions}
\label{examples}
We will apply Th. \ref{main} to the cases when $\A\subset\B$ is one of the
subnets in Lemma \ref{lr} and Lemma \ref{diagonal}, and the extension
$\A\subset \C$ corresponds to a conformal inclusion. \par
Let us start with the normal 
subnet $\A_{{SU}(2)_{10}}\subset \A_{{SU}(20)_{1}}$ as in
Lemma  \ref{lr} and extension $\A_{{SU}(2)_{10}}\subset 
\A_{{Spin}(5)_1}$. Note that the spectrum of $
\A_{{SU}(2)_{10}}\subset 
\A_{{Spin}(5)_1}$ is $L(0)+L(3)$ where $L(0)$ is the 
vacuum representation and $L(3)$ is the representation of  ${SU}(2)_{10}$ 
corresponding to the spin $3$ representation of $SU(2)$. 
By Th. \ref{main} we conclude that there is a
mirror extension $\A_{{SU}(10)_{2}}\subset \tilde C$ whose 
spectrum is $L(2\Lambda_0)+ L(\Lambda_3+\Lambda_7)$: here 
$\Lambda_3+\Lambda_7$ is the representation of ${SU}(10)_{2}$ 
corresponding to ``$(1,3)$'' in the notation of Th.  \ref{main}, and is 
determined uniquely by the branching rules for   
 $\A_{{SU}(2)_{10}}\subset \A_{{SU}(20)_{1}}$.
\par
If we choose instead  $\A_{{SU}(2)_{10}}\subset \A_{SU(2)_{K_1}}
\otimes\A_{SU(2)_{K_2}}\otimes...\otimes \A_{SU(2)_{K_l}}, K_1+...+K_l=10$ 
as the
normal subnet in the previous paragraph and apply  Th.  \ref{main}, 
we obtain mirror extensions of the coset $\tilde \A_{{SU}(2)_{10}}$
whose spectrum is determined by   Th.  \ref{main}. When $l=2,K_1=1,K_2=9$, 
the mirror extension is the extension labeled by $(A_{10}, E_6)$ in \cite{KL}. \par
Similarly we can obtain mirror extensions corresponding to $SU(2)_{28},
SU(3)_{5}, SU(3)_9$  
and $SU(3)_{21}$ cases. As in the case of $SU(2)_{10}$,
we obtain a finite number of mirror extensions. We note that in the case of
$SU(2)_{28}$ and the normal extension is 
$\A_{SU(2)_{28}}\subset  \A_{SU(2)_{27}}\otimes \A_{SU(2)_{1}}$, 
the mirror extension is the extension labeled by $(A_{28},E_8)$ in 
 \cite{KL}. \par
Next we construct infinite series of mirror extensions by using 
equations (\ref{c1}),  (\ref{c2}), (\ref{c3}). We note that if we 
choose the normal subnet as in Lemma \ref{lr}, by Th. \ref{main} 
and \S4.2 of \cite{Xu1}, we conclude that the extensions 
corresponding to ${SU}(N)_{N+2}\subset {SU}(N)_{
\frac{N(N+1)}{2}}$   
 ${SU}(N+2)_{N}\subset {SU}(N)_{
\frac{(N+2)(N+1)}{2}}$ are mirror of each other, while
the  extension corresponding to 
${SU}(N)_{N}\subset {Spin}(N)_{
N^2-1}$ is the mirror of itself.  
To obtain 
new extensions, we choose the normal subnets as in Lemma \ref{diagonal}.
We note that the spectrum of the extensions corresponding to 
the conformal inclusions of (\ref{c1}),  (\ref{c2})
and  (\ref{c3}) are given by \cite{LL} and \S4.2 of \cite{KW}. 
Applying Th. \ref{main}, we obtain three 
infinite series of mirror extensions of the diagonal cosets 
$SU(N)_{N+2}\subset SU(N)_{K_1} \times SU(N)_{K_2}\times...\times
 SU(N)_{K_l}, K_1+...K_l=N+2$, 
$SU(N)_{N-2}\subset SU(N)_{K_1'} \times SU(N)_{K_2'}\times...\times
SU(N)_{K_m'}, K_1'+...K_m'=N-2$,
$SU(N)_{N} \subset  SU(N)_{M_1} \times SU(N)_{M_2}\times...\times
 SU(N)_{M_p}, M_1+...M_p=N$,
corresponding to the conformal inclusions of (\ref{c1}),  (\ref{c2})
and  (\ref{c3}) respectively. The spectrum of these extensions are determined
as in Th. \ref{main}. \par            
Note that in all the new mirror extensions $\tilde \A\subset \tilde \C$ 
constructed in this subsetion, $\tilde \A$ corresponds to a Vertex 
Operator Algebra (VOA) (cf. \cite{FZ}) denoted by $\tilde \A_{voa}$, 
and in fact  $\tilde \A_{voa}$ correspond to
either affine Kac-Moody algebras or cosets. We will use the same notation 
as in Th. \ref{main} to label the representations of $\tilde \A_{voa}$.
Based on the close relations
between nets and VOAs as implied by \S2.2 of \cite{Xu3}, we conjecture
that each $\tilde \C$ correspond to a completely rational VOA (cf.\cite{Z})
denoted by  $\tilde \C_{voa}$ which contains $\tilde \A_{voa}$, such that
the branching rules of $\tilde \C_{voa}$ when restricted to 
$\tilde \A_{voa}$ is given by the spectrum as in Th. \ref{main}.     
 
\end{document}